\documentclass{article}
\usepackage{amssymb}
\usepackage{graphics}
\usepackage{amsthm}

\begin{document}

\newtheorem{thm}{Theorem}[section]
\newtheorem{prop}[thm]{Proposition}
\newtheorem{lem}[thm]{Lemma}
\newtheorem{cor}[thm]{Corollary}
\newtheorem{prob}{Problem}
\newtheorem{definition}{Definition}

\title{Comparing classes of finite structures}

\author{W.\ Calvert, D.\ Cummins, J.\ F.\ Knight, and S.\ Miller}

\maketitle

\section{Introduction}

In many branches of mathematics, there is work classifying a collection of
objects, up to isomorphism or other important equivalence, in terms of nice
invariants.  In descriptive set theory, there is a body of work using a notion
of ``Borel embedding'' to compare the classification problems for various
classes of structures (fields, graphs, groups, etc.) \cite{F-S},
\cite{B-K}, \cite{H}, \cite{H-K1}, \cite{H-K2}.  In this work, each class
consists of structures with the same countable universe, say $\omega$, and
with the
same language, usually finite.  For a given finite language $L$, the class
of all
$L$-structures with universe $\omega$ has a natural topological structure, and
the other classes of $L$-structures being considered are normally Borel
subclasses of these.

A \emph{Borel embedding} of one class $K$ into another class $K'$ is a Borel
function from $K$ to $K'$ that is well defined and $1-1$ on isomorphism
types.  The notation $K\leq_B K'$ indicates that there is such an embedding.
If $K\leq_B K'$, then the classification problem for $K$ reduces to that for
$K'$.  If $K'$ has nice invariants, then we may describe $\mathcal{A}\in K$, up
to isomorphism, by determining the corresponding $\mathcal{B}\in K'$ and giving
\emph{its} invariants.  If there is no nice classification for $K$, then the
same must be true for $K'$.

\bigskip
\noindent
\textbf{Example}:  Friedman and Stanley \cite{F-S} described an embedding of
the class of undirected graphs in the class of fields of characteristic $0$.
The edge relation on an undirected graph is assumed to be irreflexive.  For
an undirected graph $\mathcal{G}$, the corresponding field is obtained by
first taking an algebraically closed field of characteristic $0$, with a
transcendence base $G$, identified with the set of graph elements, and then
restricting to the subfield that is generated by the algebraic closures of the
single elements $b\in G$, and the elements $\sqrt{b_1+b_2}$, where $b_1,b_2$
are joined by an edge in $\mathcal{G}$.

\bigskip

In the present paper, our goal is to compare classes of structures using a
notion of \emph{computable} embedding.  Like the relation $\leq_B$, our
relation $\leq_c$ is a partial order on classes of structures.  We focus
mainly, but not exclusively, on classes of finite structures.  We have some
``landmark'' classes---finite prime fields, finite linear orders, finite
dimensional vector spaces over the rationals, and arbitrary linear
orders---forming a strictly increasing sequence.  If we restrict our
attention to classes of finite structures, then the class of finite linear
orders lies on top, along with the class of finite cyclic groups and
the class of finite undirected graphs.  There are many incomparable classes
below the class of finite prime fields, and between that and the class of
finite linear orders.  If we allow classes that contain infinite
structures, then the class of undirected graphs lies on top.  The
Friedman-Stanley embedding can be turned into a computable embedding, showing
that the class of fields of characteristic $0$ is also on top.  There are many
incomparable classes between the class of finite linear orders and the
class of finite dimensional vector spaces, and between this class and
the class of all linear orders.

In the remainder of the present section, we give some conventions
and definitions.  In Section~2, we discuss various natural examples of
classes.  We show that any class of finite structures can be computably
embedded in the class of finite undirected graphs, and this can be computably
embedded in the class of finite linear orders.  In Section 3, we
characterize the classes that can be computably embedded in the finite prime
fields, and those that can be computably embedded in the finite linear
orders.  In Section 4, we show that the class of finite dimensional vector
spaces over the rationals lies strictly above the class of finite linear
orders.  Using notions related to immunity, we construct families of
$2^{\aleph_0}$ incomparable classes in various intervals.  We also produce
infinite increasing chains of classes.  Finally, in Section 5, we state some
open problems.

\subsection{Conventions}

We begin with some conventions.  The structures that we consider all have a
finite relational language, and all have universe a subset of $\omega$.  The
classes that we consider all consist of structures for a single
language.  Moreover, all classes are closed under isomorphism, modulo the
restriction that each structure has universe a subset of $\omega$.  We will
sometimes identify a structure $\mathcal{A}$ with its atomic diagram
$D(\mathcal{A})$.  We will also identify finite sequences, sentences, etc.,
with their G\"{o}del numbers.  Thus, we may say that a structure is
\emph{computable}, meaning that the set of codes for sentences in
$D(\mathcal{A})$ is computable.  All finite structures are computable,
but the infinite structures that we consider may or may not be
computable.

\subsection{Basic definitions}

There are several possible notions of \emph{computable transformation} from one
class of structures to another.  The one that we have chosen is essentially
uniform enumeration reducibility.  Recall that for $A,B\subseteq\omega$, $B$ is
\emph{enumeration reducible} to $A$ if there is a computably enumerable (c.e.)
set $\Phi$ of pairs $(\alpha,b)$, where $\alpha$ is a finite subset of
$\omega$ and $b\in\omega$, such that
$$B = \{b|(\exists\alpha\subseteq A)\,(\alpha,b)\in\Phi\}\ .$$
For a given $\Phi$ and $A$, the set $B$ is unique, and we may denote it by
$\Phi(A)$.  (For more on enumeration reducibility, see Rogers \cite{Ro}.)

Here is the definition of computable transformation that we shall
use.

\begin{definition}

Let $K$ and $K'$ be classes of structures, and let $\Phi$ be a c.e.\ set of
pairs $(\alpha,\varphi)$, where $\alpha$ is a subset of the atomic
diagram of a
finite structure for the language of $K$, and $\varphi$ is an atomic
sentence, or
the negation of one, in the language of $K'$.  We say that $\Phi$ is a
\emph{computable transformation from $K$ to $K'$} if for all $\mathcal{A}\in
K$, $\Phi(D(\mathcal{A}))$ has the form $D(\mathcal{B})$, for some
$\mathcal{B}\in K'$.  We may write $\Phi(\mathcal{A}) = \mathcal{B}$
(identifying the structures with their atomic diagrams).

\end{definition}

Note that in this definition, the output $D(\mathcal{B})$ depends only on the
input $D(\mathcal{A})$, not on the order in which it is examined.

\begin{prop}
\label{prop1.1}

Let $K,K'$ be classes of structures, and let $\Phi$ be a computable
transformation from $K$ to $K'$.  If $\mathcal{A},\mathcal{A}'\in K$, where
$\mathcal{A}\subseteq\mathcal{A}'$, then
$\Phi(\mathcal{A})\subseteq\Phi(\mathcal{A}')$.

\end{prop}
\begin{proof}  Let $\mathcal{B} = \Phi(\mathcal{A})$, and let $\mathcal{B}' =
\Phi(\mathcal{A}')$.  If $\varphi\in D(\mathcal{B})$, then there is
a finite set $\alpha\subseteq D(\mathcal{A})$ such that
$(\alpha,\varphi)\in\Phi$.
Then since $\alpha\subseteq D(\mathcal{A}')$, we have $\varphi\in
D(\mathcal{B}')$.\end{proof}

It follows from Proposition \ref{prop1.1} that if $K$ contains an infinite
strictly
increasing chain of structures (increasing under the substructure relation),
then so does $K'$.  More generally, we have the following.

\begin{cor}
\label{chains}

Let $K,K'$ be classes of structures such that there is a computable
transformation from $K$ to $K'$.  If $K$ contains a strictly increasing
chain of
structures having order type $\rho$, then so does $K'$.

\end{cor}

We are interested in computable transformations that respect isomorphism,
mapping $K/_{\cong}$ into $K'/_{\cong}$ in a $1-1$ way.

\begin{definition}

Let $K,K'$ be classes of structures.

\begin{enumerate}

\item  A \emph{computable embedding} of $K$ in $K'$ is a computable
transformation $\Phi$ from $K$ to $K'$ such that for all
$\mathcal{A},\mathcal{A}'\in K$, $\mathcal{A}\cong\mathcal{A}'$ iff
$\Phi(\mathcal{A})\cong\Phi(\mathcal{A}')$.

\item  If there is a computable embedding of $K$ in $K'$, then we write
$K\leq_c K'$.

\end{enumerate}
\end{definition}

The following proposition records two obvious, but useful, facts.

\begin{prop}
\label{prop1.3}

Let $K_1,K_2, K_1',K_2'$ be classes of structures, with $K_1'\subseteq~K_1$
and $K_2'\supseteq K_2$.  If $K_1\leq_c K_2$, via $\Phi$, then
$K_1'\leq_c K_2'$, via the same $\Phi$.

\end{prop}

To illustrate what a computable embedding actually looks like, we
return to the motivating example.

\begin{prop} [Friedman and Stanley]
\label{prop1.4}

If $K$ is the class of undirected graphs, and $K'$ is the class of fields of
characteristic $0$, then $K\leq_c K'$.

\end{prop}
\noindent
\textit{Sketch of proof}.  We describe the computable embedding $\Phi$.  First,
let $\mathcal{F}$ be a  computable algebraically closed field of characteristic
$0$, with a computable sequence $(b_k)_{k\in\omega}$ of elements that are
algebraically independent.  For an undirected graph $\mathcal{G}$ (with
universe a subset of $\omega$), let $\mathcal{F}(\mathcal{G})$ be the subfield
of $\mathcal{F}$ generated by the elements that are either in the algebraic
closure of $b_k$, for some graph element $k$, or else have the form
$\sqrt{b_i+b_j}$, where $i,j$ are distinct graph elements joined by an edge.
Now, let $\Phi$ consist of the pairs $(\alpha,\varphi)$, where $\alpha$ is  the
atomic diagram of some finite undirected graph $\mathcal{G}$ and $\varphi$ is a
sentence in the atomic diagram of $\mathcal{F}(\mathcal{G})$. Clearly,
$\Phi$ is c.e.
For any $\mathcal{A}\in K$, $\Phi(\mathcal{A}) = \mathcal{F}(\mathcal{A})$.
Therefore,
$\Phi$ is a computable transformation from $K$ to $K'$.  The fact that
$\Phi$ is
$1-1$ on isomorphism types takes some effort.  (It must be shown that for
$i,j\in\mathcal{G}$, if $i,j$  are not joined by an edge, then
$\sqrt{b_i+b_j}$
is not present in $\mathcal{F}(\mathcal{G})$.)
\qed

\bigskip
\noindent
\textbf{Notation}:  We write $\mathcal{C}$ for the set of all classes of
structures satisfying our conventions, and $\mathcal{FC}$ for the restriction
to classes of finite structures.  The relation $\leq_c$ is a partial order
on $\mathcal{C}$, and as always, we get an equivalence relation $\equiv_c$,
where $K\equiv_c K'$ iff both $K\leq_c K'$ and $K'\leq_c K$.
The equivalence classes under $\equiv_c$ are called \emph{$c$-degrees}.  The
relation $\leq_c$ on $\mathcal{C}$ induces a partial order on $c$-degrees,
which we denote also by $\leq_c$.  We write $\mathbf{C}$ for the degree
structure $(\mathcal{C}/_{\equiv_c},\leq_c)$, and we write $\mathbf{FC}$ for
the restriction of $\mathbf{C}$ to the $c$-degrees that contain elements of
$\mathcal{FC}$.

\subsection{Alternative definitions}

Effective transformations between classes of structures, of one kind or
another, occur in many places in the literature.  We mention only a sample.
First, there are notions that involve \emph{interpretation}, in which the
universe and basic relations of a structure $\mathcal{A}\in K$ are defined, in
a uniform way, in the corresponding structure $\mathcal{B}\in K'$.  This
approach has been used to show that certain theories are undecidable---see,
for example, the unpublished typescript of Rabin and Scott \cite{R-S}, or the
more recent paper of Nies \cite{N}.  Sometimes, it is necessary to have the
interpretation go both ways.  This happens, for example, in the paper of
Hirschfeldt, Khoussainov, Shore, and Slinko \cite{H-K-S-S}, with results on
``computable dimension'' (the number of isomorphic members not isomorphic by a
computable function) and ``degree spectra'' (the set of possible degrees of a
relation in isomorphic copies of a computable structure).  Some of our
computable embeddings involve interpretions, but others do not.

Another kind of effective transformation, which is used in connection with
computable structures, is a partial computable function taking indices for
computable members of one class $K$ to indices for computable members of
another class $K'$.  This approach occurs, for example, in the usual
proof that the set of computable indices for computable well orderings is
$\Pi^1_1$ complete (see Rogers \cite{Ro}).  More recently, the approach is used
in \cite{C1}, \cite{C2}, in results on the complexity of the isomorphism
relation.  We deal directly with structures, not with indices, and the
infinite structures that we consider are not necessarily computable.

We state two alternative notions of computable transformation that we
tried working with, and then discarded in favor of the one in Definition 1.
In the definition below, \emph{enumeration} reducibility is replaced by
\emph{Turing} reducibility.

\bigskip
\noindent
\textbf{Definition 1$'$}:  Let $K,K'$ be classes of structures.  Then
$\Phi = \varphi_e$ (the computable operator given by oracle machine $e$) is a
\emph{computable transformation of $K$ into $K'$} if for all $\mathcal{A}\in
K$, $\varphi_e^{D(\mathcal{A})}$ is the characteristic function of
$D(\mathcal{B})$, for some $\mathcal{B}\in K'$.

\bigskip

Definition 1$'$ would be equivalent to Definition 1 if our structures all had
universe $\omega$.  However, for most structures, there will be numbers not in
the universe.  Definition 1$'$ would have us use information about these
numbers, which strikes us as not ``structural''.

\bigskip

The other alternative definition has the feature that for a given input
structure $\mathcal{A}$, the output structure $\mathcal{B}$ depends not
just on $D(\mathcal{A})$, but on the order in which we look at the
information.

\bigskip
\noindent
\textbf{Definition 1$''$}:  Let $K,K'$ be classes of structures.  An
\emph{effective transformation} of $K$ into $K'$ is a partial computable
function $f$ such that for all $\mathcal{A}\in K$, and all chains
$(\alpha_s)_{s\in\omega}$ of finite sets such that
$D(\mathcal{A}) = \cup_s\alpha_s$, there is a structure $\mathcal{B}\in K'$
such that $D(\mathcal{B})$ is the union of a corresponding chain
$(\beta_s)_{s\in\omega}$ of finite sets, where $\beta_0 = \emptyset$, and for
all $s$, $\beta_{s+1} = f(\beta_s,\alpha_s)$.

\bigskip

From Definition 1$'$, and also from Definition 1$''$, we obtain obvious
alternative versions of Definition 2, and we get further partial orders on
$\mathcal{C}$ and $\mathcal{FC}$.  Using Definition 1$''$, we would produce
transformations that respect isomorphism by guessing at the ``global''
structure.   Definition 1$''$ may be interesting from the point of view of
computability theory, especially for classes of infinite structures.  There is
a great deal of guessing at global structure in known arguments showing that
the set of indices of computable copies of various structures is $\Sigma^1_1$
complete or $\Pi^0_\alpha$-complete (see \cite{C1}, \cite{C2}).   We chose
Definition 1 as representing a more direct computable transformation of one
structure into another, based on ``local'' structure.  Proposition
\ref{prop1.1} seems intuitively right to us, and it fails for the alternative
definitions.

\subsection{Related reducibilities}

There is quite a lot of work on the Medvedev lattice and Medvedev degrees
(see \cite{Me}, \cite{Sor}, \cite{Si}).  The setting resembles ours in some
ways.  In both cases, the basic objects are classes, and the reducibility
takes members of one class to members of another in a uniform effective way.
Our reducibility relation is uniform \emph{enumeration} reducibility, while
Medvedev reducibility is uniform \emph{Turing} reducibility.  Dyment
developed a
variant of the Medvedev lattice based on enumeration reducibility (see
\cite{Sor},\cite{D}). In the Medvedev lattice, the points are classes of
\emph{functions}, while we consider classes of structures, closed under
isomorphism.  Moreover, our computable reductions are supposed to be well
defined
and $1-1$ on isomorphism types.  This makes our setting quite
different.  It will be shown in \cite{K} that $\mathcal{C}$ is not a
lattice.

\section{Examples}

Having defined the notion of a computable embedding from one class of
structures into another, we will now investigate some natural examples of
classes of structures.  If we restrict our attention to classes of finite
structures, we find that there are two distinct $c$-degrees into which almost
all natural examples of classes of finite structures fall.  One of these
$c$-degrees is made up of those classes that are computably equivalent to the
prime fields, while the other contains classes of structures that are
computably equivalent to finite linear orders (we will call these classes $PF$
and $FLO$, respectively).  We prove that these are in fact different
$c$-degrees:

\begin{prop}
\label{prop2}

$PF\lneq_c FLO$ (i.e., $PF\leq_c FLO$ and $FLO\not\leq_c PF$).

\end{prop}
\begin{proof}  To show $PF\leq_c FLO$, we construct a computable embedding
$\Phi$.  For
each $n$, let $\mathcal{L}_n$, be the the usual linear order on the first
$n$ elements
of the natural numbers.  Let $\Phi$ be the set of pairs $(\alpha,\varphi)$
such that,
for some $n$, $\alpha$ is the atomic diagram of a field of size $p_n$
(where $p_n$ is
the $n$th prime), and $\varphi\in D(\mathcal{L}_n)$.  This set of pairs is
clearly c.e.
Note that, for all
$\mathcal{A},\mathcal{A}'\in PF$, we have
$\mathcal{A}\cong\mathcal{A}'$ iff $\Phi(\mathcal{A})\cong\Phi(\mathcal{A}')$.
Therefore, $PF$ is computably embedded in $PF$.

To show that $FLO\nleq_c PF$, assume for a contradiction that $FLO\leq_c
PF$.  Say
$\Phi$ is a computable embedding.  Let $\mathcal{A}, \mathcal{A}'$ be two
nonisomorphic members of $FLO$.  Suppose $\mathcal{A}$ has fewer elements than
$\mathcal{A}'$.  Then $\mathcal{A}$ is clearly isomorphic to a substructure of
$\mathcal{A}'$.  We may suppose that $\mathcal{A}\subseteq\mathcal{A}'$.
By Proposition \ref{prop1.1},
$\Phi(\mathcal{A})\subseteq\Phi(\mathcal{A}')$.  We also
know from Definition 2 that $\Phi(\mathcal{A})\ncong\Phi(\mathcal{A}')$.
Since no
prime field is a substructure of another, nonisomorphic prime field, we
conclude that
$FLO\nleq_c PF$.
\end{proof}

We note that in the proof of Proposition \ref{prop2} above, showing
$FLO\nleq_c PF$ only used the fact that $FLO$ contains two nonisomorphic finite
linear orders. The same proof yields the following.

\begin{cor}
\label{2linorders}

If $K$ is a class containing two nonisomorphic finite linear orders, then
$K\not\leq_c PF$.

\end{cor}

We have seen that there are at least two distinct $c$-degrees in
$\mathbf{FC}$.
Most natural examples of classes of finite structures fit into one of these two
$c$-degrees, but before we discuss some of these examples, it will be
convienent
to prove that, for classes of finite structures, the $c$-degree of $FLO$ is
the maximum element of $\mathbf{FC}$.  We first prove the following lemma.

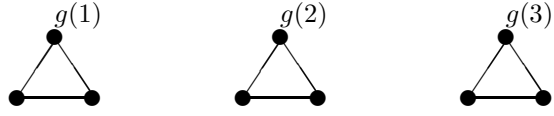
\begin{figure}
\label{fig2.1}
\setlength{\unitlength}{1mm}
\begin{picture}(60,20)(-33,0)
\linethickness{.5pt}

\put(-10,5){\circle*{2}}
\put(0,5){\circle*{2}}
\put(-5,13){\circle*{2}}
\put(-5,15){$g(1)$}

\put(-10,5){\line(2,3){5}}
\put(-10,5){\line(2,0){10}}
\put(-5,13){\line(2,-3){5}}

\put(20,5){\circle*{2}}
\put(30,5){\circle*{2}}
\put(25,13){\circle*{2}}
\put(25,15){$g(2)$}

\put(20,5){\line(2,3){5}}
\put(20,5){\line(2,0){10}}
\put(25,13){\line(2,-3){5}}

\put(50,5){\circle*{2}}
\put(60,5){\circle*{2}}
\put(55,13){\circle*{2}}
\put(55,15){$g(3)$}

\put(50,5){\line(2,3){5}}
\put(50,5){\line(2,0){10}}
\put(55,13){\line(2,-3){5}}

\end{picture}
\caption{Representing the elements $1,2,3$}
\end{figure}

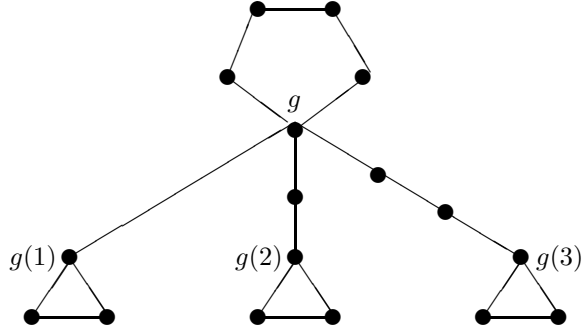
\begin{figure}
\label{fig2.3}
\setlength{\unitlength}{1mm}
\begin{picture}(60,50)(-33,0)
\linethickness{.5pt}

\put(25,30){\circle*{2}}
\put(34,37){\circle*{2}}
\put(30,46){\circle*{2}}
\put(20,46){\circle*{2}}
\put(16,37){\circle*{2}}

\put(24,33){$g$}

\put(25,30){\line(4,3){9}}
\put(30,46){\line(3,-5){5}}
\put(20,46){\line(3,0){9}}
\put(16,37){\line(2,5){4}}
\put(16,37){\line(4,-3){8}}

\put(-5,13){\line(5,3){30}}
\put(25,13){\line(0,1){10}}
\put(25,31){\line(5,-3){30}}
\put(25,21){\line(0,1){9}}

\put(-10,5){\circle*{2}}
\put(0,5){\circle*{2}}
\put(-5,13){\circle*{2}}
\put(-13,12){$g(1)$}

\put(-10,5){\line(2,3){5}}
\put(-10,5){\line(2,0){10}}
\put(-5,13){\line(2,-3){5}}

\put(25,21){\circle*{2}}

\put(36,24){\circle*{2}}
\put(45,19){\circle*{2}}

\put(20,5){\circle*{2}}
\put(30,5){\circle*{2}}
\put(25,13){\circle*{2}}
\put(17,12){$g(2)$}

\put(20,5){\line(2,3){5}}
\put(20,5){\line(2,0){10}}
\put(25,13){\line(2,-3){5}}

\put(50,5){\circle*{2}}
\put(60,5){\circle*{2}}
\put(55,13){\circle*{2}}
\put(57,12){$g(3)$}

\put(50,5){\line(2,3){5}}
\put(50,5){\line(2,0){10}}
\put(55,13){\line(2,-3){5}}

\end{picture}
\caption{Representing $R(1,2,3)$, where $R$ corresponds to $5$}
\end{figure}

\begin{lem}
\label{lem2.3}

For any class of structures $K$ in a finite relational language,
$K\leq_c~UG$, where $UG$ is the class of undirected graphs.  Moreover, if $K$
consists of finite structures, then $K\leq_c FUG$, where $FUG$ is the class of
finite undirected graphs.

\end{lem}
\begin{proof}  For each finite relational language $L$, we describe a
computable embedding $\Phi$ of the class of all $L$-structures into
$UG$.  Thus, for an arbitrary class of $L$-structures whose universes are
subsets of the natural numbers, $\Phi$ embeds the given class in $UG$.  The
embedding $\Phi$ will have the feature that if $\mathcal{A}$ is a finite
$L$-structure, then $\Phi(\mathcal{A})$ is also finite.

We begin by describing a large undirected graph $\mathcal{G}$, with finite
subgraphs
that represent possible elements of $L$-structures, and further finite
subgraphs
that represent sentences $R(a_1,\ldots,a_r)$ which may occur in the atomic
diagrams of
$L$-structures.  The graph $\mathcal{G}$ will be computable.

\bigskip
\noindent
\textbf{Subgraphs representing possible elements}

\bigskip
For each $a\in\omega$, we put into $\mathcal{G}$ a $3$-cycle $T_a$.  We
arrange
that the cycles $T_a$ are all disjoint, and we can pass effectively from
$a$ to $T_a$.
Let $g(a)$ be the least element of $T_a$.  (Figure 1 shows $T_1$, $T_2$,
and $T_3$.)

\bigskip
\noindent
\textbf{Subgraphs representing possible atomic sentences}

\bigskip
We assign to the relation symbols of $L$ distinct numbers greater than $3$.
Suppose $R$
is assigned the number $k$.  Then for each atomic sentence of the form
$\rho = R(a_1,\ldots,a_r)$, we put into $\mathcal{G}$ a subgraph
$\mathcal{G}_\rho$ consisting of a $k$-cycle together with some connecting
chains.  Say $g$ is the least element of the $k$-cycle.  We connect
$g$ to $g(a_1)$ by a chain of length $1$, adding just an edge.  We connect
$g$ to
$g(a_2)$ by a chain of length $2$, adding an intermediate point, and, in
general, we
connect $g$ to $a_i$ by a chain of length $i$, adding $i-1$ intermediate
points.  All
of the points that we have described as making up the subgraph
$\mathcal{G}_\rho$ are
distinct, and for distinct $\rho$, the subgraphs $\mathcal{G}_\rho$ are
disjoint,
except possibly for the elements $g(a)$ (in the $3$-cycles.  We arrange the
construction so that we can pass effectively from an atomic sentence $\rho$
to the
subgraph $\mathcal{G}_\rho$.  (Figure 2 gives a sample $\mathcal{G}_\rho$.)

\bigskip
The graph $\mathcal{G}$ is generated by the two families of subgraphs described
above.  For each $L$-structure $\mathcal{A}$, there is a corresponding
graph $\mathcal{G}(\mathcal{A})\subseteq\mathcal{G}$, generated by the
subgraphs $T_a$,
where $a\in\mathcal{A}$, and $\mathcal{G}_\rho$, where $\rho$ is a sentence in
$D(\mathcal{A})$ of the form $R(a_1,\ldots,a_n)$.  We note that if
$\mathcal{A}$ is
finite, then $\mathcal{G}(\mathcal{A})$ is also finite.  Now, we are ready
to define
the computable embedding $\Phi$.  This consists of the pairs
$(\alpha,\varphi)$ such
that for some finite $L$-structure $\mathcal{A}$, $\alpha = D(\mathcal{A})$ and
$\varphi\in D(\mathcal{G}(\mathcal{A}))$.

The set $\Phi$ is c.e.  For any $L$-structure $\mathcal{A}$,
$\Phi(\mathcal{A}) = \mathcal{G}(\mathcal{A})$.  It should be clear that if
$\mathcal{A}\cong\mathcal{A}'$,then $\Phi(\mathcal{A})\cong\Phi(\mathcal{A}')$. Conversely, if$\Phi(\mathcal{A})\cong\Phi(\mathcal{A}')$, then the $3$-cy
cles in
$\Phi(\mathcal{A})$ must correspond to the $3$-cycles in
$\Phi(\mathcal{A}')$.  So,
the isomorphism from $\Phi(\mathcal{A})$ onto $\Phi(\mathcal{A}')$ induces
a $1-1$
function from $\mathcal{A}$ onto $\mathcal{A}'$.  Moreover, if for some atomic
sentence $\rho = R(a_1,\ldots,a_r)$, the $3$-cycles in $\Phi(\mathcal{A})$
corresponding to $a_1,\ldots,a_r$ are attached to the subgraph
$\mathcal{G}_\rho$, indicating that $\mathcal{A}\models R(a_1,\ldots,a_r)$,
then the
corresponding $3$-cycles in $\Phi(\mathcal{A}')$ are attached to a copy of
$\mathcal{G}_\rho$, indicating that $\mathcal{A}'\models
R(f(a_1),\ldots,f(a_r))$.  It follows that $\mathcal{A}\cong\mathcal{A}'$.
Therefore,
$\Phi$ is a computable embedding of the class of all $L$-structures into
$UG$.

\end{proof}

\begin{definition}
\label{def2.1}\ \\
\begin{enumerate}
\item  A \emph{computable enumeration} of a class $K$ is a c.e.\ set
$\mathcal{E}$ of pairs $(n,\varphi)$ where $n\in\omega$ and $\varphi$ is an
atomic
sentence or the negation of one, and

\begin{enumerate}

\item for each $n$, $\{\varphi|(n, \varphi)\in\mathcal{E}\} =
D(\mathcal{A}_n)$, for
some $\mathcal{A}_n\in K$,

\item for each $\mathcal{A}\in K$, there is some $n$ such that
$\mathcal{A}_n\cong\mathcal{A}$.

\end{enumerate}
We may write $(\mathcal{A}_n)_{n\in\omega}$ instead of $\mathcal{E}$ for the
enumeration, indicating that it really is a list.

\item  An enumeration $(\mathcal{A}_n)_{n\in\omega}$ of a class $K$ is
\emph{Friedberg} if each isomorphism type in $K$ is represented just once
on the
list.

\end{enumerate}

\end{definition}

The following lemma says that $FUG$ has a computable Friedberg enumeration of a
special kind.

\begin{lem}
\label{lemFr}

There is a computable Friedberg enumeration $(\mathcal{G}_n)_{n\in\omega}$ of
$FUG$ with the feature that if $\mathcal{G}\cong\mathcal{G}_m$ and
$\mathcal{G}'\cong\mathcal{G}_n$, where $\mathcal{G}'$ is a proper extension of
$\mathcal{G}$, then $m < n$.

\end{lem}
\begin{proof}

To prove the lemma, we first define a partial
order on the class of finite undirected graphs as follows.  One graph is
greater
than another if it has more vertices than the other.  If two graphs have
the same
number of vertices, but one has more edges than the other, the one with more
edges is greater.  Two graphs are said to be \emph{equivalent} if they agree in
number of
vertices and in number of edges---equivalent graphs need not be isomorphic.
Since, for
a given set number of vertices and set number of edges, there will only be
a finite
number of different ways to arrange the edges, it is obvious that any
equivalence class
on this partial order will only have a finite number of members.  Also,
since the
number of edges a graph may contain is bounded by the number of pairs of
vertices the
graph contains, for graphs of a given number of vertices, there will only
be a finite
number of equivalence classes.

To build the Friedberg enumeration, we run through the equivalence classes, in
increasing order, and within a given equivalence class, we choose a single
representative of each isomorphism type to put into our list.  Since, for each
equivalence class, we have only a finite number of possible edge
arrangements, we
can do this effectively.  Since every finite undirected graph falls into one of
the equivalence classes, our enumeration will include every isomprphism type of
finite undirected graphs. It is a Friedberg enumeration since if two members of
$FUG$ are isomorphic, they will be in the same equivalence class, and for each
equivalence class we included just one representative of each isomorphism type.
\end{proof}

Using Lemma \ref{lemFr}, we can prove the following.

\begin{thm}
\label{thm1}

$FUG\equiv_c FLO$

\end{thm}

\begin{proof}

We must define a computable embedding $\Phi$ of $FUG$ into $FLO$.  Take the
computable Friedberg enumeration $(\mathcal{G}_n)_{n\in\omega}$ of $FUG$
with the
special feature in Lemma \ref{lemFr}.  For each $n$, let $\mathcal{L}_{n}$
be the
usual linear ordering of
$\{0,1,2... n-1\}$.  Let $\Phi$ be the set of pairs $(\alpha,\varphi)$ such
that for some $n$, $\alpha$ is the atomic diagram of a copy of
$\mathcal{G}_{n}$
and $\varphi\in D(\mathcal{L}_{n})$.  This set of pairs is clearly c.e.
For each
$\mathcal{G}\in FUG$, there is a unique $n$ such that
$\mathcal{G}\cong\mathcal{G}_n$. Then we have $\Phi(\mathcal{G}) =
\mathcal{L}_n$---here we are using the special feature of our Friedberg
enumeration, which guarantees that if $\mathcal{G}'\subseteq\mathcal{G}$, where
$\mathcal{G}'\cong\mathcal{G}_m$, then $m\leq n$.  It follows that $\Phi$
is well defined and $1-1$ on isomorphism types.
We have shown that $FUG\leq_c FLO$.  We get the fact that $FLO\leq_c FUG$
directly from Lemma \ref{lem2.3} and Proposition \ref{prop1.3}.
\end{proof}

Having shown that the $c$-degree of $FLO$ is at the top of the
$\mathbf{FC}$, we go on
to show that there are further $c$-degrees (containing classes of infinite
structures)
that lie above that of $FLO$.

\begin{thm}
\label{thm2}

The class $FVS$ of finite dimensional vector spaces over the rationals lies
strictly above $FLO$; that is, $FLO\lneq FVS$.

\end{thm}
\begin{proof}  Before proving the result, we should specify the language we
are using for vector spaces.  It is $L=\{V,F,0_{F},1_{F},+_{F},*_{F},
0_{\mathcal{V}},
+_{\mathcal{V}},*_{\mathcal{V}}\}$, where $0_{F}$,
$1_{F}$, $+_{F}$, $*_{F}$ are zero, one, addition and
mulitplication in the rationals, and $0_{\mathcal{V}}$, $+_{\mathcal{V}}$, and
$*_{\mathcal{V}}$ are the zero vector, vector addition, and multiplication of
a scalar by a vector, respectively.  Inlcuding the field symbols in our
language allows us to avoid including a separate symbol for multiplication
by each
scalar (a common approach).  Thus, our language is finite.  We make it
relational by thinking of the binary operations and constants as relations.

To prove that $FLO\leq_c FDS$, let $\mathcal{V}$ be a computable vector
space over the
rationals, with a computable sequence of basis elements $b_{1},b_{2},\ldots
$.  For
each $n$, let $\mathcal{V}_n$ be the subspace of $\mathcal{V}$ with basis
$\{b_1,\ldots,b_n\}$.  Let $\Phi$ be the set of pairs, $(\alpha,\varphi)$
such that for
some $n$, $\alpha$ is the atomic diagram of a linear ordering of size $n$ and
$\varphi\in D(\mathcal{V}_n)$.  The set $\Phi$ is clearly c.e.  Note that,
for all
$\mathcal{A},\mathcal{A}'\in FLO$, $\mathcal{A}\cong\mathcal{A}'$ iff
$\Phi(\mathcal{A})\cong\Phi(\mathcal{A}')$.  Only the number of elements in
$\mathcal{A}$ was considered in the construction of $\Phi$, so $\Phi$ will
map every
member of a given isomorphism type of $FLO$ to the same member of $FDS$.

To prove that $FDS\nleq_c FLO$, we first observe that for any finite set of
atomic
sentences $\alpha$ in the language of rational vector spaces, plus natural
numbers,
$\alpha$ is a subset of the atomic diagram of a vector space of any given
finite
dimension.  If $\alpha$ describes $n$ independent vectors in a vector space
$\mathcal{V}$, it is obvious that $\alpha$ is a subset of atomic diagrams
of vector
spaces of dimension greater than $n$, but it is also true that $\alpha$ is
a subset
of the atomic diagrams of vector spaces of dimension less than $n$.  This is
because, since $\alpha$ is finite, the sentences it contains can only
describe a
finite number of linear combinations of the $n$ vectors, saying that these
are not
$0$.  We may extend $\alpha$ to $\beta$, with a sentence saying that some
further
linear combination of two of the vectors is $0$, so that $\beta$ is a
subset of the
diagram of a vector space $\mathcal{V}'$ of dimension $n-1$.

To show that $FDS\nleq_c FLO$, assume towards a contradiction that there
exists a $\Phi$
witnessing $FDS\leq_c FLO$.  Let $\mathcal{V}$ be a two-dimensional member
of $FDS$, and
say that $\Phi(\mathcal{V})=\mathcal{L}$, where $\mathcal{L}$ is an ordering of
type $n$.  There is a finite set of pairs
$(\alpha_1,\varphi_1),\ldots,(\alpha_r,\varphi_r)$ in $\Phi$ such that
$D(\mathcal{L}) = \{\varphi_1,\ldots,\varphi_r\}$.  Then
$\alpha = \cup_{1\leq i\leq r}\alpha_i$ is a finite subset of
$D(\mathcal{V})$.  We saw
above that the set $\alpha$ is also a subset of the atomic diagram of a
rational vector
space $\mathcal{V}'$ of dimension one.  Since $\alpha\subseteq
D(\mathcal{V}')$,
$\Phi(\mathcal{V}')$ must be a linear order, say $\mathcal{L}'$, such that
$D(\mathcal{L})\subseteq D(\mathcal{L}')$.  Therefore, either
$\mathcal{L}'\cong
\mathcal{L}$, or $\mathcal{L} \subset \mathcal{L}'$. If
$\mathcal{L}'\cong\mathcal{L}$,
then $\Phi$ is not one to one on isomorphism types.  If
$\mathcal{L}\subset\mathcal{L}'$, then Proposition \ref{prop1.1} would fail for
$\Phi$.  Either way, we have our contradiction.
\end{proof}

Note that in the proof that $FDS\not\leq_c FLO$, where we used dimensions
one and two,
we could have substituted any two different dimensions.

\begin{cor}
\label{2vectorspaces}

If $K$ is a class containing vector spaces of two different dimensions, then
$K\not\leq_c FLO$.

\end{cor}

\section{General Characteristics}

The natural examples of classes described in Section 2 give rise to broader
questions regarding our ability to determine what key characteristics of those
classes are essential for their placement in our structure.  Our goal was
to give general results that determine where an arbitrary class lies in
relation to the
landmark examples, and then manipulate those results in order to construct more
examples to fill in our partial order.
For simplicity, we now will refer to the
$c$-degree containing finite prime fields as $Type\ I$ and the $c$-degree
containing
finite linear orders and finite undirected graphs as $Type\ II$.
Most generally, we know that all classes of finite structures will embed into a
class of $Type\ II$, from Lemma \ref{lem2.3}.  We would like to know what
is required
for an arbitrary class (of possibly infinite structures) to embed in a
class of $Type\
II$.  We would also like to know which classes lie above and below those of
$Type\ I$.

\subsection{Results relating to $Type\ I$}

Our examples in Section 2 suggested the abstract conditions on a class of
structures $K$
that are required for $Type\ I\leq_c K$.  The conditions involve computable
Friedberg enumerations (defined in Section 2).  It is well-known that there
are classes
with a computable enumeration but no computable Friedberg enumeration,
although we have
not been able to determine who first showed this.  There are familiar
examples, such as
the class of computable linear orderings.  Below, we construct a simple
example,
consisting of finite structures.

\begin{prop}
\label{comp f.e.}

There is a class of finite structures $K$ that has a computable
enumeration but no computable Friedberg enumeration.

\end {prop}
\begin{proof}

We want to create a class $K$ with a computable enumeration $\mathcal{E}$.  For
each $n$, the set $\{\varphi|(n,\varphi)\in\mathcal{E}\}$ should be
the diagram of some $\mathcal{A}_n\in K$, and each element of $K$ should be
isomorphic
to $\mathcal{A}_n$, for some $n$.  For each $e$, we have a requirement
$R_e$ stating that $W_e$ is not a computable Friedberg enumeration of
$K$; either $W_e$ fails to be an enumeration of $K$, or else it repeats
isomorphism types.

At each stage $s$, we have enumerated a finite subset of $\mathcal{E}$,
attempting to
take care of the first $s$ requirements.  Our strategy for $R_e$ is as
follows:   Let
$\mathcal{C}$ be an $e$-cycle, and let $\mathcal{C}^-$ be the result of
adding to
$\mathcal{C}$ an extra ``tail''---that is, an extra vertex connected to
exactly one
vertex of the cycle.  We initiate the requirement by putting into
$\mathcal{E}$ pairs $(2e,\varphi)$ for all $\varphi\in D(\mathcal{C})$ and
$(2e+1,\varphi)$ for all $\varphi\in D(\mathcal{C}^-)$.  Suppose at some
later stage
$t$, we see, for some $\mathcal{B}\cong\mathcal{C}$ and
$\mathcal{B}^-\cong\mathcal{C}^-$, and for some $m, n$,
$$\{(m,\varphi)|\varphi\in
D(\mathcal{B})\}\cup\{(n,\varphi)|\varphi\in D(\mathcal{B}^-)\}\subseteq
W_{e,t}\ .$$
Then we put into $\mathcal{E}$ any missing pairs $(2e,\varphi)$ for $\varphi\in
D(\mathcal{C}^-)$.  Thus, either $\mathcal{C}$ and $\mathcal{C}^-$ both
appear in our
enumeration, while $W_e$ fails to enumerate both, or else $\mathcal{C}^-$
appears twice
in our enumeration, and if $W_e$ is Friedberg, then in includes some
extension of
$\mathcal{C}$ not on our list.  Note that at each stage $s$ we initiate
Requirement
$R_s$, and we also look at $W_{e,s}$ to see if any requirements $R_e$, for
$e < s$,
require our adjustment.  This procedure clearly yields a class $K$ with a
computable
enumeration $\mathcal{E}$ but with no computable Friedberg enumeration.
\end{proof}

Now we have the following requirement for $Type\ I\leq_c K$.

\begin{thm}
\label{thm3.2}

For any class $K$, $Type\ I\leq_c K$ iff there is an infinite computable
Friedberg enumeration $(\mathcal{A}_n)_{n\in\omega}$ of a subclass of $K$.

\end{thm}
\begin{proof}  First, suppose $Type\ I\leq_c K$, witnessed by the compuable
embedding $\Phi$ of finite prime fields to elements of $K$.  For each
$n\in\omega$, we effectively produce a prime field of size $p_n$ (where
$p_n$ is
the $n^{th}$ prime), and we let $\mathcal{A}_n$ be the output of $\Phi$ on this
field.  Then $(\mathcal{A}_n)_{n\in\omega}$ is an infinite computable Friedberg
enumeration of a subclass of $K$.

Now, suppose that$(\mathcal{A}_n)_{n\in\omega}$ is an infinite computable
Friedberg enumeration of a subclass of $K$.  Then $Type\ I\leq_c K$, witnessed
by the computable embedding $\Phi$ that takes the prime field of size $p_n$ to
$\mathcal{A}_n$.  Clearly, $\Phi$ is well-defined and one-to-one on isomorphism
types because the enumeration of the subclass of $K$ was Friedberg.
\end{proof}

Now, we know what is required for a class $K$ to have a class of $Type\ I$
embed in it, but we would like to know what is required for $K$ to embed
into a class of $Type\ I$.  Our observation of the structure of the
representatives of the $Type\ I$ classes motivated the following definition.

\begin{definition}
\label{def2}

We say that $K$ has the \emph{substructure property} if no $\mathcal{A}_1\in K$
is isomorphic to a substructure of $\mathcal{A}_2\in K$ unless
$\mathcal{A}_1\cong\mathcal{A}_2$.

\end{definition}

\begin{prop}
\label{thm3.3}

If $K$ is a class of structures and $K\leq_c Type\ I$, then $K$ has the
substructure property.

\end{prop}
\begin{proof}  Say $\Phi$ witnesses the embedding, and suppose that we have
$\mathcal{A}_1\subseteq\mathcal{A}_2$, both in $K$, with
$\Phi(\mathcal{A}_1)=\mathcal{B}_1$ and $\Phi(\mathcal{A}_2)=\mathcal{B}_2$.
By Proposition~\ref{prop1.1}, $\mathcal{B}_1\subseteq\mathcal{B}_2$.
Since $\mathcal{A}_1\ncong\mathcal{A}_2$, we have
$\mathcal{B}_1\ncong\mathcal{B}_2$.
\end{proof}

The converse of Proposition \ref{thm3.3} does not hold.  More is needed for $K$
to embed into $Type\ I$ than just having the substructure property.  The
difficulty encountered in trying to embed a class of structures into $Type\ I$,
even with the substructure property, was that nonisomorphic structures may
still have a common substructure.  In understanding the classes of $Type\
I$, the
following definition is helpful.

\begin{definition}
\label{def3}

For $\mathcal{A}$ a structure in the language of $K$, $\mathcal{B}\in K$,
$\mathcal{A}$ is a \emph{characteristic
substructure of $\mathcal{B}$ for $K$} if and only if $\mathcal{A}$ is a
substructure of $\mathcal{B}$ and for any $\mathcal{C}\in K$ with $\mathcal{A}$
isomorphic to a substructure of $\mathcal{C}$, we have
$\mathcal{B}\cong\mathcal{C}$.  When this holds, we write
$\mathcal{A}\sqsubseteq\mathcal{B}$.

\end{definition}

The idea is that when we have seen a characteristic substructure, no
further information is needed.  With this definition, we develop the
following result.

\begin{thm}
\label{thm3.4}

Let $K$ be a class of structures.  Then the following are equivalent:

\begin{enumerate}

\item $K\leq_c Type\ I$.

\item There is a computably enumerable set $\mathcal{S}$ of pairs
$(\mathcal{A},n)$, where $\mathcal{A}$ is a finite structure in the language
of $K$, $n\in\omega$, and the following conditions hold.

\begin{enumerate}

\item For all $\mathcal{B}\in K$ there is a pair
$(\mathcal{A},n)\in \mathcal{S}$ such that $\mathcal{A}\subseteq\mathcal{B}$.

\item If $(\mathcal{A},n),(\mathcal{A}^{\prime},n^{\prime})\in\mathcal{S}$ and
$\mathcal{B}, \mathcal{B}^{\prime}\in K$, with
$\mathcal{A}\subseteq\mathcal{B}$ and
$\mathcal{A}^{\prime}\subseteq\mathcal{B}^{\prime}$, then $n=n^{\prime}$
if and only if $\mathcal{B}\cong\mathcal{B}^{\prime}$.

\end{enumerate}

\item There is a c.e.\ set $\mathcal{S}^*$ of pairs $(\varphi,n)$, where
$\varphi$ is an existential sentence in the language of $K$, $n\in\omega$, and

\begin{enumerate}

\item for each $\mathcal{B}\in K$ there exists $(\varphi, n)\in\mathcal{S}^*$
such that $\mathcal{B}\models\varphi$,

\item if $(\varphi, n), (\varphi^{\prime}, n^{\prime})\in\mathcal{S}^*$ and
$\mathcal{B}, \mathcal{B}^{\prime}\in K$, with $\mathcal{B}\models\varphi$ and
$\mathcal{B}^{\prime}\models\varphi^{\prime}$, then
$n=n^{\prime}$ if and only if $\mathcal{B}\cong\mathcal{B}^{\prime}$.

\end{enumerate}

\item There is a computable sequence $(\varphi_n)_{n\in\omega}$,
where each $\varphi_i$ is a computable $\Sigma_1$ sentence, such that

\begin{enumerate}

\item for all $\mathcal{B}\in K$ there exists $n$ such that
$\mathcal{B}\models\varphi_n$,

\item if $\mathcal{B},\mathcal{B}^{\prime}\in K$, with
$\mathcal{B}\models\varphi_n$ and
$\mathcal{B}^{\prime}\models\varphi_{n^{\prime}}$,
then $\mathcal{B}\cong\mathcal{B}^{\prime}$ if and only if $n=n^{\prime}$.

\end{enumerate}
\end{enumerate}
\end{thm}

\noindent
\textbf{Remark}:  Note that item 2 is stating that if we see
$(\mathcal{A},n)\in\mathcal{S}$ and $\mathcal{A}\subseteq\mathcal{B}$,
then $\mathcal{A}\sqsubseteq\mathcal{B}$.

\begin{proof}  First, to show that $1\Rightarrow2$, we suppose (without loss of
generality) that $K\leq_c PF$.  Let $\Phi$ witness the embedding.  We look for
$(\alpha_1,\varphi_1),...,(\alpha_k,\varphi_k)\in\Phi$ where
$\{\varphi_1,...,\varphi_k\} = D(\mathcal{F})$ for
$\mathcal{F}\cong\mathbb{F}_{p_n}$
and where there is a finite structure $\mathcal{A}$ in the language of $K$
such that
$\cup\alpha_i\subseteq D(\mathcal{A})$.  Then we put
$(\mathcal{A},n)\in\mathcal{S}$.  This $\mathcal{S}$ satisfies the desired
properties.
Next, to show $2\Rightarrow3$, we convert the given $\mathcal{S}$ into the
required
$\mathcal{S}^*$ as follows.  Whenever $(\mathcal{A},n)\in\mathcal{S}$, put
$(\varphi, n)\in\mathcal{S}^*$ where $\varphi$ is a natural existential
sentence saying
that there exist elements forming a copy of $\mathcal{A}$.

We get $3\Rightarrow 4$ immediately, letting $\varphi_n$ be the disjunction
of the
existential sentences such that $(\varphi,n)$ is in the given $S^*$.
Finally, we show $4\Rightarrow1$.  Let $(\mathcal{B}_n)_{n\in\omega}$ be a
uniformly
computable family of fields, where $\mathcal{B}_n\cong\mathbb{F}_{p_n}$.
Let $\Phi$
consist of the pairs $(\alpha(\vec c),b)$, where $\alpha(\overline{c})$ is
obtained from
a disjunct $(\exists\overline{u})\,\alpha(\overline{u})$ of $\varphi_n$, by
replacing
the tuple of variables $\overline{u}$ by a tuple of constants from
$\omega$, and
$\varphi\in D(\mathcal{B}_n)$.  Then $\Phi$ witnesses $K\leq_c Type\ I$.
\end{proof}

Motivated by the fact that many of our natural examples of classes of
structures had computable enumerations, we considered how Theorem
\ref{thm3.4} would
change if we considered only classes with this feature.  We obtained the
following
simpler result.

\begin{thm}
\label{thm3.5}

Suppose $K$ is a class of structures with a computable enumeration.  Then the
following are equivalent:

\begin{enumerate}

\item $K\leq_c Type\ I$.

\item There is a computable sequence $(\mathcal{A}_n)_{n\in\omega}$ such that
for all $n$, there exists $\mathcal{A}\in K$ such that
$\mathcal{A}_n\sqsubseteq\mathcal{A}$, and for all $\mathcal{A}\in K$, there is
a unique $n$ such that $\mathcal{A}_n\subseteq\mathcal{A}$.

\end{enumerate}

\end{thm}
\begin{proof} To show $1\Rightarrow 2$, we start with the set of pairs
$\mathcal{E}$
forming an enumeration of $K$. Say $\mathcal{E}_m$ is the structure with
$$D(\mathcal{E}_m) = \{\varphi|(m,\varphi)\in\mathcal{E}\}\ .$$
Let $\Phi$ be a computable embedding of $K$ into $PF$. For each $m$, we
look for a
finite set of pairs in $\Phi$, say
$(\alpha_1,\varphi_1),\ldots,(\alpha_k,\varphi_k)$,
such that $\mathcal{A}_n\models\alpha_k$, for $1\leq k\leq n$, and
$\{\varphi_1,\ldots,\varphi_k\} = D(\mathcal{F})$ where $\mathcal{F}$ is a
finite prime
field.  Assuming that the prime field is new; i.e., for all $k < m$,
$\Phi(\mathcal{E}_m)$ is not isomorphic to $\mathcal{F}$, we take a finite
substructure of $\mathcal{E}_n$ satisfying all $\alpha_k$, and we add this to
our list.  The sequence $(\mathcal{A}_n)_{n\in\omega}$ has the desired
properties.
To show $2\Rightarrow1$, we start with a sequence
$(\mathcal{A}_n)_{n\in\omega}$ as
in 2, and we let $S$ consist of the pairs $(\alpha,n)$, where $\alpha$ is
the atomic
diagram of a copy of $\mathcal{A}_n$.  Now, by Theorem \ref{thm3.4}, we have
$K\leq_c PF$.
\end{proof}

Together, Theorems \ref{thm3.2} and \ref{thm3.4} give us a clear picture of
the requirements for a class $K$ to have $Type\ I\leq_c K$ and $K\leq_c
Type\ I$.

\subsection{Results relating to $Type\ II$}

We would like a result saying when a class of possibly infinite structures will
embed in a class of $Type\ II$. The next result is similar to Theorem
\ref{thm3.4} in that the structures are distinguished by sentences describing
isomorphism types of substructures. To state the new result, we need the
following
definition (see \cite{A-N}, \cite{H1}, or the book \cite{A-K}).

\begin{definition}
\label{def4}

A \emph{computable $\Sigma_1$ formula} is a c.e.\ disjunction of finitary
existential forumulas, with a fixed tuple of free variables.  A
\emph{computable
$\Sigma_1$ sentence} is a computable $\Sigma_1$ formula with no free
variables.

\end{definition}

\begin{thm}
\label{thm3.6}

Let $K$ be a class of structures for the finite relational language~$L$.  Then
the following are equivalent:

\begin{enumerate}

\item $K\leq_c Type\ II$.

\item There is a c.e.\ set $\mathcal{S}$ of pairs
$(\mathcal{A},\mathcal{B})$ where $\mathcal{A}$ is a finite $L$-structure and
$\mathcal{B}$ is a finite linear ordering, such that

\begin{enumerate}

\item for any $\mathcal{C}\in K$, there exists
$(\mathcal{A},\mathcal{B})\in\mathcal{S}$ such that
$\mathcal{A}\subseteq\mathcal{C}$, and
$(\mathcal{A},\mathcal{B})$ is \emph{sufficient for $\mathcal{C}$}, in the
sense that
for $(\mathcal{A}^{\prime},\mathcal{B}^{\prime})\in\mathcal{S}$, if
$\mathcal{A}^{\prime}\subseteq\mathcal{C}$, then
$\mathcal{B}^{\prime}\subseteq\mathcal{B}$,

\item for $\mathcal{C},\mathcal{C}^{\prime}\in K, if
(\mathcal{A},\mathcal{B}),(\mathcal{A}^{\prime},\mathcal{B}^{\prime})$ are
elements of $\mathcal{S}$ sufficient for
$\mathcal{C},\mathcal{C}^{\prime}$, respectively, then
$\mathcal{C}\cong\mathcal{C}^{\prime}$ iff
$\mathcal{B}\cong\mathcal{B}^{\prime}$.

\end{enumerate}

\item There is a computable sequence $(\varphi_n)_{n\in\omega}$ of
computable $\Sigma_1$ sentences such that

\begin{enumerate}

\item for all $\mathcal{A}\in K$, there is some $n$ such that
$\mathcal{A}\models\varphi_n\ \&\ \neg\varphi_{n+1}$,

\item for all $\mathcal{A}\in K$ and all $n$,
$\mathcal{A}\models\varphi_{n+1}\rightarrow\varphi_n$,

\item for all $\mathcal{A},\mathcal{A}^{\prime}\in K$, if
$\mathcal{A}\ncong\mathcal{A}^{\prime}$, then there is some $n$ such that
$\varphi_n$ is true in only one of $\mathcal{A},\mathcal{A}^{\prime}$.

\end{enumerate}

\end{enumerate}

\end{thm}
\begin{proof}  To show $1\Rightarrow2$, suppose we have some $\Phi$ witnessing
the embedding.  We put into $\mathcal{S}$ the pairs
$(\mathcal{A},\mathcal{B})$, where $\mathcal{A}$ is a finite $L$-structure
and $\mathcal{B}$ is a finite linear ordering such that if
$D(\mathcal{B})=\{\varphi_1,...,\varphi_k\}$, there are pairs
$(\mathcal{A}_i,\varphi_i)\in\Phi$ with
$\mathcal{A}_i\subseteq\mathcal{A}$.  Now
$\mathcal{S}$ is a c.e.\ set of pairs with the properties needed for 2.

To show $2\Rightarrow3$, we take the given $\mathcal{S}$, and for each $n$,
we let
$\varphi_n$ be the disjunction, over the pairs
$(\mathcal{A},\mathcal{B})\in\mathcal{S}$
such that $\mathcal{B}$ has order type at least $n$, of existential
sentences saying
that there are elements forming a copy of $\mathcal{A}$.  This gives a
computable
sequence of computable $\Sigma_1$ sentences with the desired properties.

Finally, to show $3\Rightarrow1$, suppose that we have a computable sequence
$(\varphi_n)_{n\in\omega}$ of computable $\Sigma_1$ sentences satisfying the
three properties in 3.  Let $\mathcal{L}_n$ be the usual ordering on
$\{0,1,\ldots,n-1\}$ (as before).  Let $\Phi$ consist of the pairs
$(\alpha,\varphi)$ such that for some $n$,
\begin{enumerate}

\item $\alpha = D(\mathcal{A})$, for some finite structure $\mathcal{A}$ in the
language of $K$,

\item $\mathcal{A}\models\varphi_n$, and

\item $\varphi\in D(\mathcal{L}_n)$.

\end{enumerate}
Clearly, $\Phi$ is c.e. Moreover, we can see that $K\leq_c FLO$ via $\Phi$.
We note that if $\mathcal{A}\in K$, then $\Phi(\mathcal{A}) =
\mathcal{L}_n$, where
$n$ is greatest such that $\mathcal{A}\models\varphi_n$.
\end{proof}

Using these results, in the next section we will construct examples of
classes that fall into places in our partial order that no previous example
occupied.

\section{The structure of the partial order $\leq_c$}

In this section, we look at the partial order $\leq_c$ on $\mathcal{FC}$
(classes of finite structures) and $\mathcal{C}$ (all classes).  We begin
with $\mathcal{FC}$.  It may not have been clear, at first face, that
$\mathcal{FC}$ should have more than one $\equiv_c$-class.  However, we showed
in Proposition \ref{prop2} that there are at least two, which we called $Types\
I$ and $II$.  We are about to describe many, many more, showing that the
partial
order $(\mathcal{FC},\leq_c)$ is not only nontrivial, but highly complex.

\begin{definition}

We say that the classes $K$ and $K'$ are \emph{incomparable}, and we write
$K\perp K$, if $K\nleq_c K'$ and $K'\nleq_c K$.

\end{definition}

The result below says that there are many inequivalent classes below $Type\ I$.

\begin{prop}
\label{icpf}

There is a family of classes $(K_f)_{f\in 2^\omega}$ such that for all
$f\in~2^\omega$,
$K_f\leq_c PF$, and for
$f,g\in 2^\omega$, if $f \neq g$, then $K_f\perp K_g$.

\end{prop}
\begin{proof}  We assure that $K_f\leq_c PF$ by making $K_f\subseteq PF$ (see
Proposition \ref{prop1.3}).  There is a natural $1-1$ correspondence between
natural numbers and isomorphism types of prime fields---let the number $n$
correspond to the type of $\mathbb{F}_{p_n}$.  Then each set $A\subseteq\omega$
corresponds to the class $K_A\subseteq PF$ consisting of the fields of type
$\mathbb{F}_{p_n}$, for
$n\in A$.

To obtain a family $(K_f)_{f\in 2^\omega}$ of incomparable
subclasses of $PF$, we shall construct a family $(A_f)_{f\in\omega}$ of subsets
of $\omega$ with some special properties related to immunity.
Recall that a set $A\subseteq\omega$ is \emph{immune} if it is infinite and
has no infinite c.e.\ subset (see Soare \cite{Soa}).  We now define a
stronger property, for pairs of sets.

\begin{definition}

Let $X\subseteq\omega$.  For $A,B\subseteq\omega$, we say that $A$ and $B$ are
\emph{$X$ bi-immune} provided that for any $X$-computable function $f$ with
infinite range, there is some $a\in A$ such that $f(a)\notin B$, and
there is some $b\in B$ such that $f(b)\notin A$.  We say that $A$ and $B$ are
\emph{bi-immune} if they are $X$ bi-immune for computable~$X$.

\end{definition}

If $A$ and $B$ are $X$ bi-immune, then it is clear that neither has
any infinite $X$-computably enumerable subset, and further that no
partial $X$-computable function takes one to an infinite subset of
the other.  In a sense, $A$ and $B$ are $X$-immune with respect to each
other.  We obtain a pair of incomparable classes below $PF$ by taking a
bi-immune pair of sets $A,B$ and forming the classes $K_A,K_B$.  It is not
difficult to see that $K_A\perp K_B$.  Suppose $K_A\leq_c K_B$ via $\Phi$.
We could convert $\Phi$ into a partial computable function $f$ that maps $A$
injectively into $B$.  Let $(\mathcal{A}_n)_{n\in\omega}$ be a uniformly
computable sequence of fields such that $\mathcal{A}_n$ has type
$\mathbb{F}_{p_n}$.  Let $f(a) = b$ iff when we apply $\Phi$ to the input
$\mathcal{A}_a$, we get output describing a field of type
$\mathbb{F}_{p_b}$.

To produce the family of classes $(K_f)_{f\in 2^\omega}$ required for
Proposition \ref{icpf}, it is enough to produce a family of sets
$(A_f)_{f\in 2^\omega}$ which are pairwise bi-immune.  We prove the
following.

\begin{lem}
\label{biimmune}

For any set $X$, there exists a family $(A_f)_{f\in 2^\omega}$ such that for
any distinct $f,g\in 2^\omega$, $A_f$ and $A_g$ are $X$ bi-immune.

\end{lem}
\begin{proof}  We determine the sets $A_f$ in stages.  At stage $s$, we
associate with each $\tau\in 2^s$ a disjoint pair of finite sets
$A_\tau, A_\tau^*$, such that if $\nu\subseteq\tau$, then $A_\nu\subseteq
A_{\tau}$ and $A_\nu^*\subseteq A_{\tau}^*$.  For each $f\in 2^\omega$, we will
take $A_f$ to be the union of the sets $A_\tau$, for
$\tau\subseteq f$.  We have the following requirements.

\bigskip

\begin{tabular}{rl}

$Q_{e}$: & For all $f\in 2^\omega$, $|A_f|\geq e$\\

$R_{<e,\sigma>}$: & For all $f,g \in 2^\omega$ such that $f\supseteq\sigma
\hat{\rule{5pt}{0pt}} 0$ and $g\supseteq\sigma\hat{\rule{5pt}{0pt}}1$, if
$ran(\varphi^X_e)$ is \\
\ & infinite, then $\varphi^X_e[A_f] \nsubseteq \varphi^X_e[A_g]$, and
$\varphi^X_e[A_g] \nsubseteq \varphi^X_e[A_f]$

\end{tabular}

\bigskip
\noindent
We make a list of these requirements, with the feature that if Requirement
$s$ has
the form $R_{<e,\sigma>}$, then $|\sigma|\leq s$.  At stage $s$, we will
determine
$A_\tau$ and $A_\tau^*$ for all $\tau$ of length $s$, so as to guarantee
satisfaction of the first $s$ requirements.

We begin by letting $A_{\emptyset} = A^*_{\emptyset} = \emptyset$.  At
stage $s+1$,
we consider Requirement~$s$.  If it has the form $Q_{e}$,
then for each $\tau$ of length $s$, we take a number $k$ not in $A_\tau \cup
A_\tau^*$.  We let $A_{\tau\hat{\ }0} = A_{\tau\hat{\ }1} =   A_\tau \cup
\{k\}$,
and $A_{\tau\hat{\ }0}^* = A_{\tau\hat{\ }1}^* = A_\tau^*$.
Now, suppose Requirement $s$ has the form $R_{<e,\sigma>}$, where
$|\sigma|\leq s$.  For each $\tau$ of length $s$, we let $A_{\tau\hat{\
}0}$ and $A_{\tau\hat{\ }1}$ include the elements of $A_\tau$, and we
let $A_{\tau\hat{\ }0}^*$ and $A_{\tau\hat{\ }1}^*$ include the elements
of $A_\tau^*$.  We may add further elements as follows.  Suppose there exist
$k,k',m,m'$ such that $\varphi^X_e(k) = m$ and $\varphi^X_e(m') = k'$,
where $k\not= k'$, $m\not= m'$, and $k,k',m,m'$ are not in any of the sets
$A_\tau$, $A_\tau^*$, for $\tau$ of length $s$.  If $ran(\varphi^X_e)$ is
infinite, then there will exist such $k,k',m,m'$.  For each pair
$\tau\supseteq\sigma\hat{\ }0$, $\tau'\supseteq\sigma\hat{\ }1$ at level $s$,
we add $k$ to $A_{\tau\hat{\ }0}^*$ and $A_{\tau\hat{\ }1}^*$, and we add $m$
to $A_{\tau'\hat{\ }0}^*$ and $A_{\tau'\hat{\ }1}^*$.  Similarly, we add $m'$
to $A_{\tau'\hat{\ }0}$ and $A_{\tau'\hat{\ }1}$, and we add $k'$ to
$A_{\tau\hat{\ }0}^*,A_{\tau\hat{\ }1}^*$.

We have described the construction.  When we form the sets $A_f = \cup_{\sigma
\subseteq f} A_\sigma$, as planned, each of the requirements is satisfied,
and the
conclusion of the lemma holds.  We note that the construction could be
carried out
using a $\Delta^0_2(X)$ oracle, so that the assignment of finite sets
$A_{\tau}$ and
$A_{\tau}^*$ to $\tau\in 2^{<\omega}$ is $\Delta^0_2(X)$.
\end{proof}

Having completed the proof of Lemma \ref{biimmune}, we have also completed
the proof of Proposition \ref{icpf}.
\end{proof}

There are also incomparable degrees that are not below $Type\ I$.  These may be
obtained by using $\Delta^0_2$ bi-immune sets and letting them determine linear
orders instead of prime fields.

\begin{prop}
\label{iclo}

There is a family of classes $(K_f)_{f\in 2^\omega}$ such that
for all $f\in~2^\omega$, $K_f\leq_c FLO$ and $K_f\perp PF$, and for distinct
$f,g\in 2^\omega$, $K_f\perp K_g$.

\end{prop}
\begin{proof}  We show that any pair of $\Delta^0_2$ bi-immune sets gives rise
to a pair of subclasses of $LO$ that are incomparable with each other and
with $PF$.
Then to obtain the family of  classes $(K_f)_{f\in 2^\omega}$ with the required
properties, we apply Lemma \ref{biimmune} to get a family $(A_f)_{f\in
2^\omega}$ of
pairwise $\Delta^0_2$ bi-immune sets, and let $K_f$ be the set of linear
orders whose
sizes are in $A_f$.

Let $A$ and $B$ be $\Delta^0_2$ bi-immune sets.  Let $K_1$ be the class of
linear orders
whose sizes are members of $A$, and let $K_2$ be the class of linear orders
whose
sizes are members of $B$.  By Proposition \ref{prop1.3}, $K_i\leq_c FLO$.  We
must show that $K_1\perp K_2$.  Suppose not, say $K_1\leq_c K_2$ via $\Phi$.
We convert $\Phi$ into a partial $\Delta^0_2$ function $f$ that maps $A$
injectively into $B$.  Let $(\mathcal{L}_n)_{n\in\omega}$ be a uniformly
computable
family of orderings, where $\mathcal{L}_n$ has type $n$.  We let $f(a) = b$ if
$\Phi$ takes $\mathcal{L}_a$ to an ordering of type $b$.  Using $\Delta^0_2$,
we can determine, for each input structure $\mathcal{A}_a$, the full atomic
diagram of the output structure $\Phi(\mathcal{L}_a)$.  Now, $f$ maps $A$
injectively into $B$, contradicting the assumption that $A,B$ are
$\Delta^0_2$ bi-immune.

We must show that $K_i\perp PF$.  The fact that $K_i\nleq_c PF$ follows from
Corollary \ref{2linorders}.  Suppose $PF\leq_c K_1$ via $\Phi$.  Let
$(\mathcal{A}_n)_{n\in\omega}$ be a uniformly computable sequence of
fields, where
$\mathcal{A}_n$ has type $\mathbb{F}_{p_n}$.  Then we have an injective
$\Delta^0_2$ function $g$ from $\omega$ into $A$, defined so that $g(n)$ is
the number
of elements in $\Phi(\mathcal{A}_n)$.  This contradicts the immunity
assumptions.
Therefore, $PF\nleq_c K_i$.
\end{proof}

We still have not shown that there are classes properly between $Types\ I$
and $II$.  On one hand, it seems that the only difference between these two
types is whether, in building the structure, we can tell whether
we're done or not. Thus, it would not be surprising to see that there was
simply nothing in between.  On the other hand, there is a sense in which
$Type\ I$ looks analogous to a computable degree, and $Type\ II$ to a complete
c.e.\ degree, so it is also reasonable to think that there are things between
them. It turns out that this second argument may be closer to the truth.

\begin{prop}
\label{prop4.4}

There is a pairwise incomparable family of classses $(\hat{K}_f)_{f\in
2^\omega}$ such that $PF\lneq_c \hat{K}\lneq_c FUG$, where
$FUG$ is the class of finite undirected graphs.

\end{prop}
\begin{proof}  For simplicity, the discussion here will show only how to
produce a
single class $\hat{K}$.  The construction of $2^{\aleph_0}$ incomparable
classes $\hat{K}_f$ would follow the outline of Lemma \ref{biimmune}.  The
class
$\hat{K}$ will be made up of finite graphs.  This
guarantees that $\hat{K}\subseteq FUG$.  Let
$(\mathcal{C}_n)_{n\in\omega}$ be a uniformly computuble sequence of cyclic
graphs of
size $n$.  To guarantee that $PF\leq_c\hat{K}$, we
include the graphs of isomorphism type $\mathcal{C}_{2n}$, for all
$n\in\omega$.
This suffices, since we have a computable embedding that takes fields of type
$\mathbb{F}_{p_n}$ to $\mathcal{C}_{2n}$.

To guarantee that $\hat{K}\nleq_c PF$, we satisfy the following requirements.

\bigskip

\begin{tabular}{rl}

$R_e$: & $W_e$ does not witness that $\hat{K}\leq_c PF$.

\end{tabular}

\bigskip

The strategy for $R_e$ is as follows.  We give $W_e = \Phi$ input
$\mathcal{C}_{2e+1}$,
and see if it produces a prime field as output (we could determine this using a
$\Delta^0_2$ oracle).  If not, then we put all copies of $\mathcal{C}_{2e+1}$
into $\hat{K}$.  If, given input $C_{2e+1}$, $\Phi$ produces as
output some finite prime field, then we do not put copies of
$\mathcal{C}_{2e+1}$ into $\hat{K}$.  Instead, we add all copies of two
different extensions of $\mathcal{C}_{2e+1}$.  We could take these to be the
result of adding a single new vertex, and either connecting it to one of
the vertices
of $\mathcal{C}_{2e+1}$, or not.

We have described all elements of the class $\hat{K}$.  We have satisfied
each requirement $R_e$---either $C_{2e+1}\in\hat{K}$, and $\Phi$ does not map
it to a finite prime field, or else $\hat{K}$ contains two nonisomorphic
extensions of $C_{2e+1}$ (neither isomorphic to a substructure of the other),
and $\Phi$ fails to map them to nonisomorphic prime fields, since
the diagram of $\Phi(C_{2e+1})$ is contained in the output for both
extensions.  Finally, we must show that $FUG\nleq_c\hat{K}$.  For this, we use
Corollary~\ref{chains}, noting that there are arbitrarily large finite
increasing
chains of graphs, and there are no chains of structures in $\hat{K}$ of length
greater than one.\end{proof}

In fact, there is an infinite chain of classes between $Types\ I$ and $II$, all
incomparable with $\hat{K}$.  These, and more examples to come, are formed by
starting with a class of $Type\ I$ (for aesthetic reasons, usually cyclic
graphs)
and adjoining subclasses of a $Type\ II$ class.  For instance, we can build
all the examples constructed here using only cyclic graphs and chains (simply
connected graphs in which each vertex is connected to at most two others).  The
following will be useful.  Clearly any class containing only these elements
is a subclass of the set of finite graphs, and is thus reducible to $Type\ II$.

\begin{prop}
\label{prop4.5}

Let $\hat{K}$ be as in Proposition \ref{prop4.4}.  Then there is a sequence of
classes $(K_n)_n\in\omega$ such that
$$PF\equiv_c K_0\lneq_c K_1\lneq_c K_2\lneq_c\ldots \leq_c FLO\ ,$$
and for all $n >0$, $K_n \perp \hat{K}$.

\end{prop}
\begin{proof}

For the moment, we let $K_n$ consist of all finite prime fields
and all chains of length $j$ for $j\leq 2n$.  This will make it easier for us
to refer back to the construction of $\hat{K}$ in the proof of Proposition
\ref{prop4.4}.  At the end of the proof, we shall replace the prime fields by
cyclic graphs.

Since $PF = K_0$, we have $PF\equiv_c K_0$.  For all $n$, we have
$K_n\leq_c K_{n+1}$, since $K_n\subseteq K_{n+1}$ (using Proposition
\ref{prop1.3}).  We must show that $K_{n+1}\nleq_c K_n$.  If
$K_{n+1}\leq_c K_{n}$, witnessed by $\Phi$, then $\Phi$ maps at least two
nonisomorphic chains to finite prime fields.  We may suppose that one of these
chains is a substructure of the other.  Since no prime field is a substructure
of another, Corollary \ref{chains} gives a contradiction.  Thus, there
is no such $\Phi$, and $K_n\lneq_c K_{n+1}$.  We have $K_n\leq_c FLO$, for all
$n$, just because all of the structures in $K_n$ are finite---in
Section 2, we saw that all classes of finite structures can be computably
embedded in $FLO$.

We must show that for all $n > 0$, $K_n\perp\hat{K}$.  We
get the fact that $K_n\nleq_c\hat{K}$ using Corollary \ref{chains}.  The
class $K_n$ contains a chain of structures of length at least $2$, while
$\hat{K}$ has no such chains.  Finally, we show that $\hat{K}\nleq
K_n$.  Suppose $\hat{K}\leq_c K_n$, witnessed by $\Phi = W_e$.  We constructed
$\hat{K}$ so that either $\mathcal{C}_{2e+1}$ is in $\hat{K}$ and $\Phi$ fails
to map $\mathcal{C}_{2e+1}$ to a finite prime field, or else $\hat{K}$ contains
two nonisomorphic extensions of $\mathcal{C}_{2e+1}$, while $\Phi$ gives
output for
both that contains the diagram of the same finite prime field (so $\Phi$ cannot
map these extensions to either finite prime fields or chains).

Now, $\Phi$ cannot map $\mathcal{C}_{2e+1}$ to a finite prime field, and
must instead
map it to one of the chains in $K_n$.  Recall that there are infinitely
many different
indices $e'$ for the same c.e.\ set $\Phi$.  If $\Phi = W_e'$, then the
argument above
shows that $\Phi$ must map $\mathcal{C}_{2e'+1}$ to one of the chains in
$K_n$.  Since
there are only finitely many isomorphism types of chains in $K_n$, and
$\Phi$ must map
infinitely many nonisomorphic cyclic graphs to them, $\Phi$ cannot be $1-1$
on isomorphism types.  Thus, $\hat{K}\nleq_c K_n$.

At this point, we replace the finite prime fields in each class $K_n$ by the
finite cyclic graphs.  The resulting class is computably equivalent to $K_n$,
but it satisfies the convention that all structures in a class have the same
finite relational language.  We justify the temporary violation of our
convention that all members of a class have a common language by noting that
we could have have substituted finite cyclic graphs for the prime fields in
the construction in Proposition \ref{prop4.4}, and in the classes
$K_n$.\end{proof}

Now, any class consisting of all finite cyclic graphs and infinitely many
finite chains will lie strictly above $K_n$, for all $n$, and will be bounded
above by $Type\ II$.  We will use exactly this sort of class to produce
many more
incomparable classes between the chain of $K_n$'s and $Type\ II$.

\begin{prop}
\label{prop4.6}

Let $(K_n)_{n\in\omega}$ be as in Proposition \ref{prop4.5}.  Then there is
a family of pairwise incomparable classes $(H_f)_{f\in2^\omega}$ lying above
all $K_n$ and below $Type\ II$.

\end{prop}
\begin{proof}  We show how to produce two classes $H,H'$.  Each class will
contain all finite cyclic graphs.  We add finite chains so as to satisfy
the following requirements:

\bigskip

\begin{tabular}{rl}

$R_e$: & $W_e$ does not witness $H\leq_c H'$\\

$R_e'$:  & $W_e$ does not witness $H'\leq_c H$

\end{tabular}

\bigskip
\noindent
We make a list of the requirements and satisfy them in order.  The
strategy for $R_e$ is as follows.  Let $\Phi = W_e$.  For earlier
requirements, we will have decided, for finitely many
$k$, whether or not to put chains of length $k$ into
$H,H'$.  Take $n$ greater than any of these $k$.  Then $n$ is also an upper
bound on
the number of chains already in $H'$. We add chains to $H$ so that there is an
increasing sequence
$\mathcal{L}_{0}\subseteq\ldots\subseteq\mathcal{L}_{n}$ of length $n+1$.  If
$\Phi$ does not map these $\mathcal{L}_{i}$ to an increasing sequence of
chains, then the requirement is already satisfied.  (If $\Phi$ maps some
$\mathcal{L}_i$ to a prime field, then we would have a contradiction of
Corollary~\ref{chains}.)  If $\Phi$ maps the $\mathcal{L}_i$ to an increasing
sequence of chains, then at least one, say the chain of length $m$, is not
already in $H'$.  We satisfy the requirement by keeping chains of length $m$
out of $H'$.  The strategy for $R_e'$ is the same, so it is clear that we
can produce two incomparable classes $H,H'$, lying above all $K_n$ and below
$Type\ II$.

In satisfying any one requirement, we make only finitely many decisions about
which chains do and do not belong to a given class.  Therefore, we could
use the
same strategy to produce a family $(K_f)_{f\in 2^\omega}$ of incomparable
classes.  We would follow the outline in the proof of Lemma 4.2.
\end{proof}

The results so far, for classes of finite structures, are summarized in
Figure~3.  Finite orders lie at the top, along with finite undirected
graphs. Finite cyclic groups, and finite simple groups lie there too.  Prime
fields lie strictly lower.  The empty class obviously lies on the bottom.  The
classes consisting of copies of a single finite structure lie just above
that---equivalent to classes consisting of copies of a single computable
structure.  The numbers 4.3, 4.6, 4.1 refer to propositions showing the
existence of large incomparable families, and 4.5 refers to the proposition
producing a chain.  The question marks indicate places where there may or may
not be a class, lying below certain classes and above certain others.

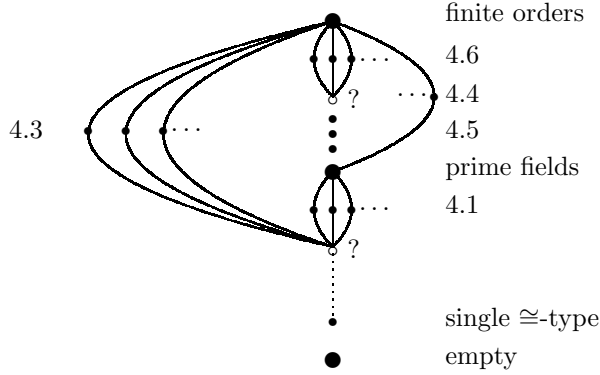
\begin{figure}
\label{fig1}
\setlength{\unitlength}{1mm}
\begin{picture}(60,50)(-10,0)
\linethickness{.5pt}

\put(45,45){\circle*{2}}
\put(60,45){finite orders}

\put(42.5,40){\circle*{1}}
\put(45,40){\circle*{1}}
\put(47.5,40){\circle*{1}}
\put(48.5,39.8){$\dots$}

\put(60,39.5){$4.6$}

\put(45,34.5){\circle{1}}
\put(47,33.5){?}

\put(53.5,35.05){\dots}
\put(58.5,35){\circle*{1}}
\put(60,34.25){$4.4$}
\qbezier(45,25)(72,35)(45,45)

\qbezier(45,35)(40,40)(45,45)
\qbezier(45,35)(45,40)(45,45)
\qbezier(45,35)(50,40)(45,45)

\put(45,32){\circle*{1}}
\put(45,30){\circle*{1}}
\put(45,28){\circle*{1}}

\put(60,29.5){$4.5$}

\put(2,29.5){$4.3$}

\put(45,25){\circle*{2}}
\put(60,24.5){prime fields}

\put(42.5,20){\circle*{1}}
\put(45,20){\circle*{1}}
\put(47.5,20){\circle*{1}}
\put(48.5,20){\dots}

\qbezier(45,15)(40,20)(45,25)
\qbezier(45,15)(45,20)(45,25)
\qbezier(45,15)(50,20)(45,25)

\put(60,19.5){$4.1$}

\put(45,14.5){\circle{1}}
\put(47,13.5){?}

\qbezier(45,15)(-20,30)(45,45)
\qbezier(45,15)(-10,30)(45,45)
\qbezier(45,15)(0,30)(45,45)

\put(12.5,30.5){\circle*{1}}
\put(17.5,30.5){\circle*{1}}
\put(22.5,30.5){\circle*{1}}
\put(23.5,30.5){\dots}

\qbezier[8](45,6)(45,10)(45,14)

\put(45,5){\circle*{1}}

\put(45,0){\circle*{2}}

\put(60,4.5){single $\cong$-type}
\put(60,-.5){empty}

\end{picture}
\caption{Classes of finite structures}
\end{figure}

Of course, the structure is dramatically enriched when we also consider
classes containing infinite structures.  We have seen that the class $FVS$ of
finite dimensional vector spaces over the rationals sits strictly above
$Type\ II$.  First, let us show that there are classes below $FVS$ that are not
below $Type\ II$.

\begin{prop}
\label{prop4.7}

There is a family of classes $(K_f)_{f\in 2^\omega}$,
such that for all\linebreak $f\in 2^\omega$, $K_f\leq_c FVS$ and $K_f$ is
incomparable with $Types\ I$ and $II$, and for distinct $f,g\in 2^\omega$,
$K_f\perp K_g$.

\end{prop}
\begin{proof}  Let $A,B\subseteq\omega$ be a pair of $\Delta^0_3$ bi-immune
sets, and
let $K, K'$ be the classes of $\mathbb{Q}$-vector spaces whose dimensions
are members of
$A$, $B$, respectively.  We have $K,K'\leq_c FVS$, by
Proposition~\ref{prop1.3}.
Next, we show that $K,K'$ are incomparable with $Types\ I$ and $II$.  By
Corollary
\ref{2vectorspaces}, no class containing vector spaces of two different
dimensions
embeds in $FLO$.  Therefore $K,K'\nleq_c FLO$, and it follows that
$K,K'\nleq_c PF$.  Let $(\mathcal{A}_n)_{n\in\omega}$ be a uniformly
computable sequence of fields, where $\mathcal{A}_n\cong\mathbb{F}_{p_n}$.
Suppose $PF\leq_c K$ via $\Phi$.  From $\Phi$, we will obtain a $\Delta^0_3$
function $f$ mapping $\omega$ injectively into $A$.  We let $f(n)$ be the
dimension of $\Phi(\mathcal{A}_n)$.  Note that the relation
$\{(n,m)|dim(\Phi(\mathcal{A}_n))\geq m\}$ is $\Sigma^0_2$.  From this, it
is clear
that $f$ is $\Delta^0_3$.  This contradicts the immunity properties of $A$.
Therefore, $PF\nleq_c K$, and similarly $PF\nleq_c K'$.  It follows that
$FLO\nleq_c K,K'$.

Similarly, we show that $K\perp K'$.  Let $(\mathcal{V}_n)_{n\in\omega}$ be a
uniformly computable family of vector spaces, where $\mathcal{V}_n$ has
dimension $n$.  If $K\leq_c K'$ via $\Phi$, then we would have a partial
$\Delta^0_3$ function $f$ mapping $A$ injectively into $B$.  We let $f(a)$
be the
dimension of $\Phi(\mathcal{V}_a)$---assuming that this is a finite-dimensional
vector space, we can apply a $\Delta^0_3$ procedure to find the dimension.

Now, Lemma \ref{biimmune} yields a family $(A_f)_{f\in 2^\omega}$ of pairwise
$\Delta^0_3$ bi-immune sets.  For each $f\in 2^\omega$, we let $K_f$ be the
class of vector spaces of dimension in $A_f$.  The argument above shows that
this family has all of the properties needed in Proposition \ref{prop4.7}.
\end{proof}

The interval between $Type\ II$ and $FVS$ admits further complexity. There
is an
$\omega$-chain with an infinite antichain above it, all in this interval.

\begin{lem}
\label{lem4.8}

If $K = LO\cup FVS$, then $K\leq_c FVS$.

\end{lem}
\begin{proof}  We define a computable embedding $\Phi$ that takes a linear
order of size $n$ to a vector space of dimension $2n+1$, and takes a vector
space of dimension $n$ to one of dimension $2n$.  Only the last part requires
verification.  We partition $\omega$ into three infinite computable sets
$A,B,C$, and we let $f,g,h$ be injective computable functions mapping
$\omega$ into
$A,B,C$, respectively.  Suppose $\alpha$ is a finite set of atomic sentences
and negations of atomic sentences (appropriate to be included in the
diagram of a vector space over $\mathbb{Q}$).  Say
$\alpha$ describes distinct vectors $n_0,\ldots,n_k$, where $n_0$ is the zero
vector.  We modify $g$, letting $g(n_0) = f(n_0)$.  Let $\alpha^*$ consist
of sentences $\varphi(f(n_0),\ldots,f(n_i))$ and
$\varphi(g(n_0),g(n_1),\ldots,g(n_k))$, where
$\varphi(n_0,\ldots,n_k)\in\alpha$, plus further sentences saying
$f(n_i) + g(n_j) = h(<n_i,n_j>)$, for $i,j\not= 0$, and sentences generated
from these by the axioms for vector spaces.  For example, if we have the
sentence $q\cdot n_i = n_j$, where $i\not= 0$ and $q$ is a non-zero rational,
then we obtain the sentence $q\cdot h(<n_i,n_i>) = h(<n_j,h_j>)$.  We put into
$\Phi$ the pairs $(\alpha,\varphi)$, where $\alpha$ is as described, and
$\varphi$
is in the corresponding set $\alpha^*$.
\end{proof}

\begin{prop}
\label{prop4.9}

There is a sequence of classes $(J_n)_{n \in \omega}$ such that
\[LO \leq_c J_0\lneq_c J_1\lneq_c J_2\ldots \leq_c FVS\ .\]

\end{prop}
\begin{proof} We define $J_n$ to be the class containing all finite linear
orders and all rational vector spaces of dimension at most $2n$.  If
we wish, we can consider these as structures in a single language, one that
enables us to distinguish a linear order from a vector space.  Clearly,
$LO = J_0$.  The proof that $FVS \nleq_c LO$ (from Theorem \ref{thm2}) shows
that $J_1\nleq_c LO$.  By Proposition \ref{prop1.3}, we have $J_n \leq_c
J_{n+1}$.  Finally, if $J_{n+1}\leq_c J_n$ via $\Phi$, then $\Phi$ would map
at least two isomorphism classes of vector spaces to
isomorphism classes of linear orders, which is again impossible.
\end{proof}

In Proposition \ref{prop4.6}, we obtained incomparable classes strictly between
$Types\ I$ and $II$ by adding to the class all of finite prime fields the
linear
orders of selected sizes.  Below, we use the same idea to obtain incomparable
classes above all $J_n$ and below $FVS$.

\begin{prop}
\label{prop4.10}

Let $(J_n)_{n\in\omega}$ be as in Proposition \ref{prop4.9}.
There exist pairwise incomparable classes $(G_f)_{f\in 2^omega}$ such that
for all
$n\in\omega$ and all $f\in 2^\omega$, we have $J_n\leq_c G_f\leq_c FVS$.  (Then
$J_n\lneq_c G_f\lneq_c FVS$.)

\end{prop}

\begin{proof}
Each class
$G_n$ will contain all finite linear orders.  We add vector spaces of selected
dimensions so as to satisfy the following requirements.

\bigskip
\begin{tabular}{rl}

$R_{\langle e,\sigma\rangle}$: & for all
$f\supseteq\sigma\hat{\ }0,g\supseteq\sigma\hat{\ }1$,
$W_e$ does not witness
$G_f \leq_c G_g$.\medskip\\

$R'_{\langle e,\sigma\rangle}$: & for all $f\supseteq\sigma\hat{\ }0,
g\supseteq\sigma{\ }1$, $W_e$ does not witness $G_g\leq_c G_f$.

\end{tabular}

\bigskip
\noindent
We have a list of all requirements.  At each stage $s$, we have decided, for
finitely many pairs $(i,n)$, whether to put vector spaces of dimension $n$ in
$G_f$, for $f\supseteq\tau$, for $\tau$ of length $s$.  We write $G_\tau$
for the
class reflecting the decisions made up to stage $s$.  In the end, we will let
$G_f$ be the union of the sets $G_\tau$, for $\tau\subseteq f$.

We satisfy the requirements  in order.  At stage $s+1$, we consider
Requirement~$s$.
Say this is $R_{\langle e,\sigma\rangle}$, and let $\Phi = W_e$.  Take $n$
greater than
the dimension of any vector spaces considered so far.  For all
$\tau\supseteq\sigma\hat{\ }0$ of length $s+1$, we put into $G_\tau$ the vector
spaces of dimension $n+i$ for $i\leq n$.  If $\Phi$ does not map one
of the newly added vector spaces to anything in $FLO\cup FVS$, then the
requirement
is trivially satisfied. If $\Phi$ embeds $G_\tau$ into $FLO\cup FVS$, then
by the
argument in Theorem \ref{thm2} or Corollary \ref{2vectorspaces}, since
$G_\tau$ contains
vector spaces of at least two different finite dimensions, $\Phi$ cannot
map any vector
space in $G_\tau$ to a finite linear order.  Thus, $\Phi$ must map one of
the $n+1$ new
vector spaces to a vector space of dimension greater than $n$.  We satisfy the
requirement by keeping the vector spaces of this dimension out of $G_\nu$,
for all
$\nu\supseteq\sigma\hat{\ }1$ of length $s+1$.
\end{proof}

\begin{figure}
\label{fig2}
\setlength{\unitlength}{1mm}
\begin{picture}(60,49)(-15,5)
\linethickness{.5pt}

\put(45,55){\circle*{2}}
\put(55,54.5){undirected graphs}

\put(43.5,51){$\equiv_c$?}

\put(45,49){\circle*{2}}
\put(55,48.5){linear orders}

\put(36,44){4.12}

\put(45,40){\circle*{2}}
\put(55,39.5){vector spaces}

\qbezier(45,40)(45,44)(45,48)

\put(42.5,35){\circle*{1}}
\put(45,35){\circle*{1}}
\put(47.5,35){\circle*{1}}
\put(48.5,35){$...$}

\put(33.5,33){$4.10$}

\qbezier(45,30)(40,35)(45,40)
\qbezier(45,30)(45,35)(45,40)
\qbezier(45,30)(50,35)(45,40)

\put(45,29.5){\circle{1}}
\put(47,28){$?$}

\put(45,27){\circle*{1}}
\put(45,25){\circle*{1}}
\put(45,23){\circle*{1}}

\put(37,24.5){$4.9$}

\put(45,20){\circle*{2}}
\put(55,19.5){finite orders}

\qbezier[20](50.5,26.5)(56,37)(45,49)

\put(47,22){\circle*{1}}
\put(48.5,23.5){\circle*{1}}
\put(50,25){\circle*{1}}

\put(55,30){$4.11$}

\put(45,15){\circle*{2}}
\put(55,14.5){prime fields}

\put(5,15){\circle*{1}}
\put(10,15){\circle*{1}}
\put(15,15){\circle*{1}}
\put(16,15){$\ldots$}

\put(0,20){$4.7$}

\qbezier(5,15)(18,40)(45,40)
\qbezier(10,15)(20,40)(45,40)
\qbezier(15,15)(22,40)(45,40)

\put(45,10){\circle*{2}}
\put(55,9.5){empty}

\end{picture}
\caption{Classes of structures, possibly infinite}
\end{figure}
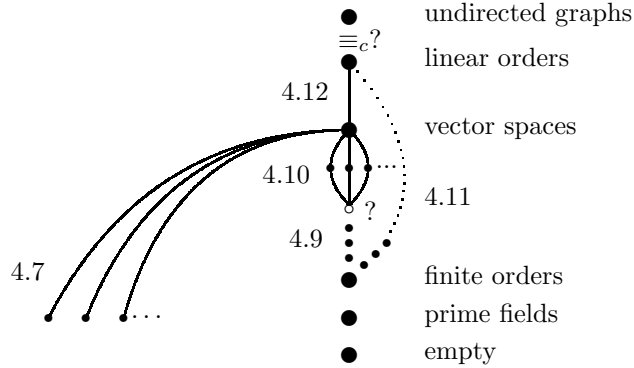

There are other classes not below $FVS$.  Of course, any class is embeddable
in the class of infinite graphs.  Also, we have an increasing sequence of
classes that is not below $FVS$.  For each ordinal $\alpha < \omega_1$, let
$LO^\alpha$ denote the set of well orders of order type less than
$\alpha$.

\begin{prop}
\label{lonemb}

If $\alpha < \beta < \omega_1$, then $LO^\alpha \lneq_c LO^\beta$.
Also, for $\omega < \alpha$, $LO^\alpha \lneq_c FVS$.

\end{prop}
\begin{proof}

If $\alpha <\beta$, then Proposition \ref{prop1.3} yields the
fact that $LO^\alpha \leq_c LO^\beta$.  There are representatives of the
order types in $LO^\alpha$ forming a chain of structures of length $\alpha$.
Then by Corollary \ref{chains}, if $\alpha < \beta$, $L^\beta\nleq L^\alpha$,
and if $\alpha > \omega$, then $L^\alpha\nleq FVS$.
\end{proof}

This result also suffices to show that the class of infinite linear orders
does not lie below the finite dimensional vector spaces.  However, we have the
following:

\begin{prop}
\label{prop4.12}

If $LO$ is the class of linear orders (possibly infinite), and $FVS$ is the
class
of finite dimensional $\mathbb{Q}$-vector spaces, then $FVS\leq_c LO$.

\end{prop}
\begin{proof}

Each vector space $\mathcal{V}$ will correspond to a substructure of
$\omega\cdot\eta$, in which for each $n\leq dim(\mathcal{V})$, we have densely
many copies of the finite linear order of size $n$, and also densely many
copies of $\omega$.  Clearly if we can describe a computable transformation
that behaves this way, it will be well-defined and injective on isomorphism
types.

Let $\mathcal{B}$ be the lexicographic ordering on $\mathbb{Q}\times\omega$.
We partition $\mathbb{Q}$ computably into dense subsets
$\mathbb{Q}_{\vec{a}}$, corresponding to finite sequences $\vec{a}$ of natural
numbers.  Given a finite set $\alpha$ of atomic sentences and
negations of atomic sentences describing an $n$-tuple of vectors
$\vec{a}$, let $\alpha^*$ describe the restriction of
$\mathcal{B}$ to $\mathbb{Q}_{\vec{a}}\times\omega$, if $\alpha$ contains
evidence that $\vec{a}$ is dependent, and $\mathbb{Q}_{\vec{a}}\times n$,
otherwise.  Let $\Phi$ consist of the pairs $(\alpha,\varphi)$, where
$\varphi\in\alpha^*$.
\end{proof}

Thus, the situation when we include classes containing infinite
structures looks something like Figure 4.  Undirected graphs lie on top.
Linear orders may or may not be equivalent to undirected graphs.  Finite
dimensional vector spaces over $\mathbb{Q}$ lie strictly below linear
orders, and finite linear orders lie strictly below vector spaces.  The
numbers 4.7, 4.10, etc., indicate the propositions being
illustrated.

\section{Problems}

In this section, we list some open problems.

\begin{prob}
\label{prob1}

Is the class of graphs computably equivalent to the class of linear orderings?

\end{prob}
The class of graphs (including infinite as well as finite ones) lies at the
top of our partial ordering. Problem \ref{prob1} asks whether there is a
computable embedding of graphs in linear orderings.

\begin{prob}
\label{prob2}

Is there a ``natural'' class $K$, consisting of finite structures of infinitely
many different isomorphism types, such that $K$ is not computably equivalent to
either finite prime fields or finite linear orderings?  In particular, is
there a natural example of a class properly between these two?

\end{prob}

We have results characterizing those classes that computably embed in the
finite prime fields, and also in the finite linear orderings.  For classes
that embed in the finite dimensional vector spaces over $\mathbb{Q}$, we can
give some necessary conditions, but we have no characterization.

\begin{prob}
\label{prob3}

Characterize the classes $K$ such that $K\leq_c FVS$.

\end{prob}

We have not entirely sorted out the differences between the definition of
$\leq_c$ that we chose and the two alternative definitions.  We can show that
the partial ordering obtained from Definition 1$''$ differs from
$\leq_c$.  We do not know about the partial ordering obtained from
Definition 1$'$.

\begin{prob}
\label{prob5}

Is it true that for any classes of structures $K,K'$, $K\leq_c K'$ iff there
is a computable operator $\Phi = \varphi_e$ of the kind in Definition 1$'$,
taking
$\mathcal{A}\in K$ to $\mathcal{B}\in K'$ such that
$\varphi_e^{D(\mathcal{A})} =
\chi_{D(\mathcal{B})}$, in a way that is well-defined and $1-1$ on isomorphism
types?

\end{prob}

\end{document}